\documentclass{amsart}

\usepackage[T1]{fontenc}
\usepackage[utf8]{inputenc}

\theoremstyle{plain}

\usepackage[all]{xy}
\usepackage{times}
\usepackage{verbatim,bm}
\usepackage{enumitem}
\usepackage[english]{babel}
\usepackage{amssymb}
\usepackage{amsthm}
\usepackage{amsmath}
\usepackage{amstext}
\usepackage{multirow}
\usepackage{hyperref}
\usepackage{graphicx}
\usepackage{epstopdf}
\usepackage{mathrsfs}
\usepackage{amscd}
\usepackage{color}

\newtheorem{thm}{Theorem}[section]
\newtheorem*{thm*}{Theorem}

\newtheorem{conj}[thm]{Conjecture}

\theoremstyle{remark}

\author{Gibran Espejo-Ramos} 
\address{Instituto de Matem\'aticas, Universidad Nacional Aut\'onoma de 
M\'exico 	(UNAM), 	\'Area de la Investigaci\'on Cient\'ifica, Circuito 
exterior, Ciudad 	Universitaria, 04510, Ciudad de M\'exico, M\'exico. }
\email{gibranr@im.unam.mx}

\author{Frank Loray} 
\address{Univ Rennes, CNRS, IRMAR, UMR 6625, F-35000 Rennes, France. }
\email{frank.loray@univ-rennes.fr}

\author{Laura Ortiz-Bobadilla } 
\address{Instituto de Matem\'aticas, Universidad Nacional Aut\'onoma de 
M\'exico 	(UNAM), 	\'Area de la Investigaci\'on Cient\'ifica, Circuito 
exterior, Ciudad 	Universitaria, 04510, Ciudad de M\'exico.}
\email{laura@im.unam.mx}

\thanks{F. Loray is supported by CNRS, ANR-11-LABX-0020-0 Labex program Centre Henri Lebesgue and ANR-25-CE40-1360 project IsoMoDyn. G. Espejo-Ramos et L. Ortiz-Bobadilla are supported by Papiit (Dgapa UNAM) IN103123 and Secihti CBF-2025-I-2101.}
	
 
\begin{document}
\title[Neighborhoods of elliptic curves and  the Riemann-Hilbert correspondence]{Neighborhoods of singular elliptic curves and the Riemann-Hilbert correspondence at infinity.}
\date{\today}
\begin{abstract}
In his thesis work, the first author studies neighborhoods of singular elliptic curves
of type $I_0^*$ (Kodaira's classification) and provides the analytic classification. 
This is related via a double covering construction to the analytic invariants found by Touzet, Voronin and the second author
 for the analytic classification of neighborhoods of smooth elliptic curves.
After recalling these constructions, we explain in this note how this construction is conjecturally related to the Riemann-Hilbert
correspondence at infinity for Painlevé VI case when considering neighborhood of Okamoto divisor.
\end{abstract}

\maketitle
\tableofcontents

Let $X$ be the elliptic curve defined over $\mathbb C$ by the Legendre affine equation $y^2=x(x-1)(x-t)$, and marked point 
at infinity. We also denote by $\tilde X\simeq\mathbb C\ni z$ its universal cover, with lattice $\mathbb Z+\tau\mathbb Z$.
We will consider the elliptic involution 
$$\iota\ :\ X\to X\ ;\ (x,y)\mapsto(x,-y)$$
induced by $z\mapsto -z$ in the universal cover,
and denote by $\pi:X\to X/\iota$ the quotient map. 
This quotient can be viewed either as an orbifold $X/\iota$, or simply as $\mathbb P^1\ni x$ with marked points 
$0,1,t,\infty$.

In this note, we will consider different moduli spaces attached to $X$, and see how the elliptic involution
is acting on them. This provides the simplest cases of the Riemann-Hilbert correspondence.
Then we will consider infinite dimensional deformations of these objects, namely
germs of neighborhoods of smooth/singular elliptic curves and their classifying maps
introduced in \cite{LTV} and \cite{TheseGibran}. Finally, we focus on a four-dimensional family 
given by Okamoto surfaces, related to the Painlevé VI equation, and settle a conjecture
about the Riemann-Hilbert correspondance at infinity. We end with some open questions.

\section{An example of Atiyah-Serre}\label{sec:AtiyahSerre}

A rank-one connection on $X$ is a pair $(L,\nabla)$ where $L$ is a line bundle over $X$
and $\nabla:L\to L\otimes\Omega^1_X$ is a $\mathbb C$-linear map of sheaves 
satisfying the Leibniz rule $\nabla(f\cdot s)=s\otimes df+f\nabla(s)$ for any local function $f$
and section $s$.
Because $X$ has dimension one, $(L,\nabla)$ is automatically flat and can be locally analytically trivialized 
as the horizontal connection $d$ on holomorphic functions.  Any two local trivializations differ by a constant $\mathbb C$-linear map,
i.e. an element of $\mathrm{GL}_1(\mathbb C)\simeq\mathbb C^*$;
this gives rise to a local system on $X$, 
is an element of the cohomology 
group $H^1(X,\mathbb C^*)$. 
This is the easiest example of a moduli space of connections. 
From the exact sequence of sheaves
$$0\to\mathbb  C^*\to \mathcal O_X^*\to \Omega_X^1\to 0$$
(where $\mathcal O_X^*\to \Omega_X^1$ is given by $f\mapsto\frac{df}{f}$),
we can derive the associate long exact sequence of cohomology groups and extract the following
\begin{equation}\label{Eq:ExactSeqRk1}
0\to H^0(X,\Omega_X^1) \stackrel{\alpha}{\to} H^1(X,\mathbb C^*) \stackrel{\beta}{\to} \mathrm{Jac}(X)\to 0
\end{equation}
making $H^1(X,\mathbb C^*)$ into an affine $\mathbb C$-bundle over the curve $\mathrm{Jac}(X)\simeq X$.
One may compactify it by adding a section at infinity in order to get
a $\mathbb P^1$-bundle $\beta:P\to X$: it is one of the two undecomposable $\mathbb P^1$-bundles found 
by Atiyah in \cite{Atiyah}. The section at infinity $\sigma:X\to P$ is the unique section whose image $C:=\sigma(X)\simeq X$
has zero self-intersection in the total space $S:=\mathrm{Tot}(P)$ of the bundle. 
The group $H^1(X,\mathbb C^*)$ identifies, as a complex manifold, with the complement $S\setminus C$.
This is usually called deRham moduli space. Because we consider rank one connections, the corresponding moduli space
$H^1(X,\mathbb C^*)$ is a group, where multiplication corresponds to the tensor product $\otimes$ of line bundles
and connections. 

Viewing $X$ as a Riemann surface, we get the Riemann-Hilbert correspondence 
$$RH\ :\ H^1(X,\mathbb C^*)\longrightarrow \mathrm{Hom}(\pi_1(X),\mathbb C^*)$$
which associates to a local system its monodromy representation (again we identify $\mathrm{GL}_1(\mathbb C)\simeq\mathbb C^*$). Once we have choosen
a system of generators $(1,\tau)$ for the fundamental group, we get an identification between
the space of representations and the group $\mathbb C^*\times\mathbb C^*$.
This is usually called Betti moduli space.
We thus get a complex analytic group isomorphism
\begin{equation}\label{Eq:RHaffineRk1}
RH\ :\ H^1(X,\mathbb C^*)\longrightarrow \mathbb C^*\times\mathbb C^*.
\end{equation}
However, the two underlying varieties are not isomorphic from the algebraic point of view:
 while $\mathbb C^*\times\mathbb C^*$
is an affine variety, there is no non constant regular functions on $S$.
Indeed, such a function would extend as a rational function on $S$ with polar divisor $n\cdot C$;
one easily check that the zero divisor of such function should define a complete curve not intersecting $C$.
But there is no complete curve in $S$ since it is analytically isomorphic to 
an affine variety. This famous example is due to Serre, and provides, when varrying the analytic structure
of the curve $X$, a  one-parameter family of different projective compactifications of $\mathbb C^*\times\mathbb C^*$
(see \cite{Ueda1}).

Let us briefly mention two other facts. The maximal compact subgroup given by unitary connections/representation
is an analytic real torus providing a real section of $P\to X$, and turns to be real algebraic on representation-side,
namely $\mathbb S^1\times\mathbb S^1\subset\mathbb C^*\times\mathbb C^*$. The fibration $P\to X$
is sent to the foliation defined by $ u\partial_u+\tau v\partial_v$ where $(1:\tau)$ is commensurable
to periods of $X$.  

One can make it more explicit. The Atiyah-Serre surface $S$ can be viewed, analytically, as the quotient of 
$\mathbb C\times\mathbb P^1\ni(z,w)$ by the action of the transformation group generated by 
\begin{equation}\label{eq:ActionSerre}
\phi_1(z,w)=(z+1,w)\ \ \ \text{and}\ \ \ \phi_\tau(z,w)=(z+\tau, w-1)
\end{equation}
To view this, observe first that
any flat line bundle admits a connection with trivial monodromy $u=1$ along the first generator:
by the exact sequence (\ref{Eq:ExactSeqRk1}), it suffices to add a convenient holomorphic
one-form to a given connection in order to normalize its monodromy.
Any flat line bundle can therefore be described as the quotient $L_z$
of $\mathbb C\times\mathbb C\ni(\zeta,\lambda)$ by the action of the transformation group generated by 
\begin{equation}
\varphi_1(\zeta,\lambda)=(\zeta+1,\lambda)\ \ \ \text{and}\ \ \ \varphi_\tau(\zeta,\lambda)=(\zeta+\tau, e^{2i\pi z}\lambda).
\end{equation}
The other connections on this bundle are given by the quotient $\nabla_{z,w}$ of $\lambda\mapsto d\lambda+(2i\pi wd\zeta)\lambda$ for $w\in\mathbb C$.
After integration, we deduce that the monodromy is given by
\begin{equation}
RH\ :\ (z,w)\mapsto(u,v)=(e^{2i\pi w},e^{2i\pi(z+\tau w)}).
\end{equation}
In order to derive the moduli space $S\setminus C$, we have to identify those $(L_z,\nabla_{z,w})$ giving rise 
to the same monodromy representation, i.e. quotienting by the maps $\phi_1,\phi_\tau$.

The advantage of this formulation is that we can easily deduce the topology at infinity of Riemann-Hilbert map. 
Let us denote by $V^0:=\mathbb P^1\times\mathbb P^1$ the compactification of Betti moduli space,
and by $D=L_1\cup L_2\cup L_3\cup L_4$ the compactifying divisor, where 
$$L_1=\{u=0\},\ \ \ L_2=\{v=0\},\ \ \ L_3=\{u=\infty\},\ \ \ L_4=\{v=\infty\}$$
so that $V^0\setminus D=\mathbb C^*\times\mathbb C^*$. 
If we denote by $V_i^0$ a small tubular neighborhood of $L_i$, 
then we can check that $U_i:=RH^{-1}(V_i^0\setminus D)\subset S\setminus C$ is a transversely
sectorial open set along $C$. For instance, if $V_1^0$
$$V_1^0=\{\vert u\vert<<0\}\ \ \ \Rightarrow\ \ \ U_1=\{\Im(w)>>0\}$$
which is a sector of opening $-\pi<\arg(\frac{1}{w})<0$.
Therefore, the two loops corresponding to the periods $(1,\tau)$ are sent to 
the loops around $L_1$ and $L_2$ in $V^0\setminus D$ and, near infinity,
the loop around $C$ in $S\setminus C$ is sent to a loop turning along the cycle of curves $D$.

There is a natural non degenerate $2$-form $\omega$ on $S\setminus C$ constructed from the 
exact sequence (\ref{Eq:ExactSeqRk1}) by choosing non zero $1$-forms on the base $X$ and on the fiber;
it is well-defined up to a scalar constant. One can define it in uniformizing coordinates $(z,w)$ as $dz\wedge dw$.
It extends meromorphically on $S$ with polar divisor $2C$. Isomorphism (\ref{Eq:RHaffineRk1}) sends it to 
the $2$-form $\eta:=\frac{du}{u}\wedge\frac{dv}{v}$ (up to a scalar).

\section{Deformation of the neighborhood of the elliptic curve}

One can focus on the analytic structure of the neighborhood at infinity of the surface $S$, i.e. 
a neighborhood $U^0$ of the embedded curve $C$. Since $\phi_1,\phi_\tau$ are tangent to the identity along $w=\infty$,
we deduce that the normal bundle $N_C$ of $C=\{w=\infty\}$ is trivial. However, the neighborhood
is not linearizable: the trivial elliptic fibration on $N_C$ does not exist of $S$
because $S\setminus C$ contains not compact curve. The infinitesimal obstruction of existence
of fibrations (or more generally linear connections) on neighborhoods of compact cuves has been 
considered by Ueda in \cite{Ueda2}. In the case $N_C$ is torsion, the zero section is the (possibly multiple) fiber 
of a unique elliptic fibration on the total space $\mathrm{Tot}(N_C)$; in this case, the {\bf Ueda type} of the neighborhood
$k\in\mathbb Z_{>0}\cup\{\infty\}$ is the maximal integer for which the elliptic fibration 
extends to the $k^{th}$ infinitesimal neighborhood $J^k(U,C)$ (see also \cite{LTT}).
In the Atiyah-Serre example, one can check that Ueda type is $1$.

In \cite{LTT,LTV}, the formal and analytic classification of two-dimensional neighborhoods of elliptic curves
with torsion normal bundle is established. We say that two neighborhood germs $(U,C)$ and $(U',C)$
are conjugated if, and only if, we have a commutative diagram
\begin{equation}\label{eq:equivNeighbId}
\xymatrix{\relax
    C \ar[r]^\iota \ar[d]_{\text{id}}  & U \ar[d]^\phi \\
    C \ar[r]^{\iota'} & U'
}\end{equation}
where $\phi$ is a complex analytic diffeomorphism between appropriate representatives $U,U'$ of the germs
(we hope that the reader will not be confused by our many abuse of notations).
Each formal class contains an infinite-dimensional deformation
of different analytic classes. Let us explain the case where the Ueda type is $1$.
We call a {\bf square-type} neighborhood a two-dimensional germ of neighborhood $(V,D)$ of a (simple normal crossing)
singular curve $D$ composed of a cycle of $4$ smooth rational curves, $D=L_1\cup L_2\cup L_3\cup L_4$,
with vanishing self-intersection $L_i\cdot L_i=0$, likely as the compactifying divisor $D=\mathbb P^1\times\mathbb P^1\setminus\mathbb C^*\times\mathbb C^*$.

\begin{thm}\label{thm:EquivalenceSquare}
There is a one-to-one correspondence between analytic classes of neighborhoods $(U,C)$ 
of the elliptic curve $C$ with Ueda type $1$, and analytic classes of square-type neighborhoods $(V,D)$.
\end{thm}

Analytic classification of square-type neighborhoods $(V,D)$ can easily be described as follows.
By Kodaira-Spencer, the neighborhood germ $(V,L_i)$ of each line $L_i$ is trivial, 
analytically equivalent to the corresponding germ of neighborhood in $V^0$:
$$\psi_i\ :\ (V,L_i)\stackrel{\sim}{\to}(V^0,L_i).$$
Therefore, the neighborhood germ $(V,D)$ can be constructed by patching together
trivial neighborhood germs $(V^0,L_i)$ by means of complex analytic diffeomorphisms
$$\varphi_{i,i+1}\in\mathrm{Diff}(V^0,D,p_{i,i+1}),\ \ \ i\in\mathbb Z/4$$
fixing $p_{i,i+1}=L_i\cap L_{i+1}$ and preserving $D$ locally.  
On the other hand, these  gluing maps $\varphi_{i,i+1}$ depend on the choice 
of the trivializations $\psi_i$, that can be postcomposed by a complex analytic diffeomorphisms
$$\varphi_i\in\mathrm{Diff}(V^0,D,L_i),\ \ \ i\in\mathbb Z/4$$
preserving $D$ along $L_i$. The moduli space of analytic classes is finally described by the collection of $4$-uples $(\varphi_{i,i+1})$
up to action of $(\varphi_i)$.

For instance, $\varphi_{1,2}(u,v)=(u\cdot  a(u,v),v\cdot  b(u,v))$, where $a,b\in\mathbb C\{X,Y\}^\times$ (non vanishing),
can be replaced by $\varphi_1\circ\varphi_{1,2}\circ\varphi_2$, where 
$\varphi_1(u,v)=(u\cdot  a_1(v), v\cdot  b_1(v))$ and $\varphi_2(u,v)=(u\cdot a_2(u), v\cdot  b_2(u))$ 
with $a_1,b_1,a_2,b_2\in\mathbb C\{X\}^\times$. 

It turns out that the proof of Theorem \ref{thm:EquivalenceSquare} goes through a deformation of the 
Riemann-Hilbert map described in the previous section. 

\begin{thm} Let $(U,C)$ be a germ of neighborhood of the elliptic curve $C$ with Ueda type $1$,
and $(V,D)$ be the square type neighborhood germ associated to it by Theorem \ref{thm:EquivalenceSquare}.
Then, there is a natural complex analytic isomorphism 
$$LTV\ :\ U\setminus C\stackrel{\sim}{\longrightarrow} V\setminus D$$
between representatives $U,V$ which is asymptotic to the Riemann-Hilbert map when $(U,C)$ is 
formally equivalent to Atiyah-Serre example.
\end{thm}

The main tool in the proof is a sectorial normalization. Given a basis $(1,\tau)$ of the lattice of the elliptic curve $C$,
there is a natural decomposition $U\setminus C=U_1\cup U_2\cup U_3\cup U_4$ by transversely sectorial open sets,
together with complex analytic isomorphisms $\Psi_i:U_i\stackrel{\sim}{\to} U_I^0$ onto similar neighborhoods for Atiyah-Serre surface $S$. 
The $LTV$ map is then constructed on sectors by considering $LTV_i:=RH\circ \Psi_i\ :\ U_i\to V^0\setminus D$.
Then $LTV_i=\varphi_{ij}\circ LTV_j$, and $LTV$ is well-defined after patching $LTV_i$'s images into $(V,D)$ by means of $\varphi_{ij}$'s.

\section{Quotient by the elliptic involution}\label{sec:Quotient}

The involution $\iota$ acts on rank-one local systems considered in section \ref{sec:AtiyahSerre},
and therefore on the deRham moduli space  $S\setminus C$ as 
$$\iota\ :\ S\to S\ ;\ (z,w)\mapsto(-z,-w).$$
It extends trivially on the compact ruled surface $\beta:S\to X$ and has $8$ fixed points:
on each fiber $F_0,F_1,F_t,F_\infty$, where $F_x:=\beta^{-1}(x)$, $\iota$ has one fixed point at infinity
$p_x=F_x\cap C$, and one other fixed point $p_x'\in F_x\setminus p_x$ which corresponds to a $2$-torsion 
point of $H^1(X,\mathbb C^*)$.

The quotient surface $\beta:S/\iota\to X/\iota$ has $8$ conic points of type $A_1$. We denote
by $Bl:\Sigma\to S/\iota$ the desingularization given by blowing-up the $8$ points, and by $E_x,E_x'$ 
the corresponding exceptional divisors. Equivalently, we can first blow-up $\hat{S}\to S$ at $p_x,p_x'$
and then take the quotient: we get the same surface $\Sigma=\hat{S}/\iota$. 
The curve $C$ is invariant by $\iota$ and defines a rational curve $\Delta\subset\Sigma$.
The divisor $\Gamma:=2\Delta+ E_0+ E_1+ E_t+ E_\infty$ is a singular elliptic curve of type $I_0^*$
(Kodaira's classification) and the polar divisor of the two-form induced by $\omega$ on the quotient.
Here, $I_0^*$-type means that 
$$\Delta\cdot \Delta=-2,\ \ \ C\cdot E_x=1,\ \ \ E_x\cdot E_x=-2,\ \ \ \text{and}\ \ \ E_{x}\cdot E_{x'}=0 \ \ \ \text{if}\ x\not=x'$$
for any $x,x'\in\{0,1,t,\infty\}$.
The open set $\Sigma\setminus\Gamma$ can be interpreted 
as the moduli space of rank $2$ meromorphic connections defined as pushforward $\pi_*(L,\nabla)$ 
of the rank $1$ connections (see section \ref{sec:AtiyahSerre}). It corresponds to the Okamoto space
of initial conditions for Picard-Painlevé equation.

\begin{equation}\label{eq:quotientSerre}
\xymatrix{\relax
    {\color{white}W=}\hat{S} \ar[r]^{\hskip0.3cm\text{blow-up}} \ar@<10pt>[d]_{\hat{\Pi}}  & S \ar[d]^{\Pi} \\
   \Sigma= \hat{S}/\iota \ar[r]^{\text{\hskip0.3cmblow-up}} & S/\iota
}\end{equation}

On the other hand, the involution acts on the fundamental group, its representations and the Betti moduli space
as 
$$\iota_B\ :\ (u,v)\mapsto\left(\frac{1}{u},\frac{1}{v}\right).$$
The quotient $(\mathbb C^*\times\mathbb C^*)/\iota_B$ can be identified with the affine surface 
$$\mathcal{X}=\{(X,Y,Z)\in\mathbb C^3\ ;\ X^2+Y^2+Z^2+XYZ=4\},$$
known as Cayley cubic, and we get the explicit quotient map
$$\mathbb C^*\times\mathbb C^*\mapsto \mathcal{X}\ :\ 
(u,v)\mapsto (X,Y,Z)=\left(-u-\frac{1}{u},-v-\frac{1}{v},-uv-\frac{1}{uv}\right).$$
This surface can be interpreted as rank $2$ representations of the orbifold fundamental group
$$\rho:\pi_1^{\text{orb}}(X/\iota)\to\mathrm{SL}_2(\mathbb C)$$
that split as decomposable representation after covering $X\to X/\iota$.
From this point of view, it is a particular case of character varieties 
attached to $\mathrm{SL}_2(\mathbb C)$-representations of the fundamental group of the $4$-punctured sphere
(see \cite{BG}).

The Riemann-Hilbert map commutes with the actions of $\iota$ on both sides, 
giving rise to a complex analytic symplectic isomorphism
$\underline{RH}\ :\ \Sigma\setminus\Gamma\stackrel{\sim}{\to}\mathcal{X}$
which a particular case of the  Riemann-Hilbert correspondance for Painlevé equations,
namely Picard-Hitchin case, and the general case will be described in section \ref{sec:Painleve}.

The usual compactification $\hat{\mathcal{X}}$ of $\mathcal{X}$, by viewing $\mathcal{X}\subset\mathbb C_{(X,Y,Z)}\subset\mathbb P^3$, 
is smooth at infinity, and consists in adding a cycle of $3$ rational curves with self-intersection $-1$,
defined by $XYZ=0$ at infinity (see \cite{CL}). But we will rather work with the surface 
$\overline{\mathcal{X}}=(\mathbb P^1\times\mathbb P^1)/\iota$ 
which is obtained by blowing-down the $(-1)$-curve $Z=0$ in $\hat{\mathcal{X}}$. 
Precisely, the action of $\iota_B$ on $D$ is without fixed point, permuting
$$L_1\stackrel{\iota_B}{\leftrightarrow} L_3\ \ \ \text{and}\ \ \ L_2\stackrel{\iota_B}{\leftrightarrow} L_4$$
and the quotient is a cycle $\underline{D}=\underline{L_1}\cup\underline{L_2}$ of $2$ smooth rational curves
with $\underline{L_1}\cdot\underline{L_2}=2$ and $\underline{L_i}\cdot\underline{L_i}=0$.
We will call {\bf half-square type} this kind of divisor.
We can resume the previous discussion by the following commutative diagram:
\begin{equation}\label{eq:quotientRHAtiyahSerre}
\xymatrix{\relax
 \hat{S} \ar[r]^{\text{blow-up}} \ar[rd]_{(2:1)} &S\ar@{.>}[d]&   S\setminus C \ar@{_{(}->}[l] \ar[r]^{RH}_\sim \ar[d]_{(2:1)}  & \mathbb C^*\times\mathbb C^* \ar[d]^{(2:1)} \ar@{^{(}->}[r] & \mathbb P^1\times\mathbb P^1 \ar[d]^{(2:1)} \\
 &\Sigma&  \Sigma\setminus \Gamma \ar@{_{(}->}[l] \ar[r]^\sim_{\underline{RH}} & \mathcal{X} \ar@{^{(}->}[r] & \overline{\mathcal{X}}
}\end{equation}
The two-form $\frac{du}{u}\wedge\frac{dv}{v}$ on $C^*\times\mathbb C^*$ induces the two-form $-\frac{dX\wedge dY}{2Z+XY}$
on $\mathcal{X}$ which has simple poles along $\underline{D}$.

\section{Deformation of the quotient}\label{sec:TheseGibran}

In his PhD, the first author generalized this construction for general neighborhoods of elliptic curves $(U,C)$
that admit an elliptic involution, i.e. extending the elliptic involution $\iota$ to the neighborhood
with only $4$ isolated points, lying on $C$.

\begin{thm}[G. Espejo Ramos, \cite{TheseGibran}] 
There is a one-to-one correspondence between:
\begin{enumerate}
\item analytic classes of neighborhood germs $(U,C)$ of smooth elliptic curve
with Ueda-type $1$, together with an elliptic involution $\iota:(U,C)\to(U,C)$,
\item analytic classes of neighborhood germs $(W,\Gamma)$ 
of singular elliptic curves $\Gamma$ of Kodaira type $I_0^*$,
\item analytic classes of neighborhood germs of half-square type $(\underline{V},\underline{D})$.
\end{enumerate}
Moreover, these neighborhoods germs are connected by natural maps generalizing the diagram (\ref{eq:quotientRHAtiyahSerre})
\begin{equation}\label{eq:quotientLTV}
\xymatrix{\relax
 \hat{U} \ar[r]^{\text{blow-up}} \ar[rd]_{(2:1)} & U\ar@{.>}[d]&   U\setminus C \ar@{_{(}->}[l] \ar[r]^{LTV}_\sim \ar[d]_{(2:1)}  & V\setminus D \ar[d]^{(2:1)} \ar@{^{(}->}[r] & V \ar[d]^{(2:1)} \\
 &W&  W\setminus \Gamma \ar@{_{(}->}[l] \ar[r]^\sim_{\underline{LTV}} & \underline{V}\setminus\underline{D} \ar@{^{(}->}[r] & \underline{V}
}\end{equation}
\end{thm}
All neighborhoods are smooth neighborhoods in the statement. The left-hand-side covering is ramifying along 
the branches $E_0,E_1,E_t,E_\infty$ of the $I_0^*$-divisor $\Gamma$, while the three other ones are \'etale.

\section{The Painlevé VI family of Riemann-Hilbert maps at infinity}\label{sec:Painleve}

In the series of paper \cite{DM,Iwasaki,IIS1,IIS2}, the Riemann-Hilbert correspondance 
is described in many details for the deRham and Betti moduli spaces associated to the Painlevé VI
family of differential equations. Let us recall briefly the construction.

Let us fix $t\in\mathbb P^1\setminus\{0,1,\infty\}$ and $\boldsymbol{\kappa}=(\kappa_0,\kappa_1,\kappa_t,\kappa_\infty)\in\mathbb C$. 
One can first associate
the deRham moduli space $M_{dR}^{\boldsymbol{\kappa}}(t)$ of rank $2$ logarithmic connections
$(E,\nabla)$ where 
\begin{itemize}
\item $E\to \mathbb P^1$ is a rank $2$ vector bundle on $\mathbb P^1$,
\item $\nabla:E\to E\otimes\Omega^1(T)$ is a logarithmic connection with polar divisor $T=0+1+t+\infty$,
\item the residue $\mathrm{Res}_x\nabla$ has eigenvalues $\pm\kappa_x$ for $x=0,1,t,\infty$.
\end{itemize}
For generic $\boldsymbol{\kappa}$, i.e. if we avoid equalities of the form $\kappa_x=0$ or $\sum_x\pm \kappa_x\in\mathbb Z$,
then the moduli space of isomorphism classes $M_{dR}^{\boldsymbol{\kappa}}(t)$ of such connections has a natural
structure of complex irreducible and smooth quasi-projective variety of dimension $2$. For general parameters, 
we can obtain a similar result by considering moduli space of semi-stable parabolic connections with respect 
to some choice of weights. The moduli space is therefore given by GIT quotient, and it turns out that, in the generic
case, it depends neither on the parabolic structure, nor on the choice of weights. This construction is done 
in \cite{IIS2} for general $\boldsymbol{\kappa}$. 

It has been known since the works of R. Fuchs (1907) that Painlevé VI equation governs the isomonodromic 
deformations of rank $2$ linear differential equations (or connections) on the Riemann sphere with $4$ simple poles.
Painlevé Property tells us that Painlevé equations define non linear local systems. 
K. Okamoto proved this property in \cite{Okamoto} and explicitely described the fiber of this local system,
namely the space of initial conditions, when we fix the position of poles, i.e. $t$ in the previous setting.
This space naturally identifies with the deRham moduli space $M_{dR}^{\boldsymbol{\kappa}}(t)$.
One way to describe it is to consider the diagonal $\Delta_0\in\mathbb P^1\times\mathbb P^1$
and first blow-up the $4$ points $p_x=\mathrm{pr}_1^{-1}(x)\cap\Delta_0$, $x=0,1,t,\infty$, where 
$\mathrm{pr}_1$ denotes the first projection to $\mathbb P^1$.
Denote by $\Sigma_0\to\mathbb P^1\times\mathbb P^1$ this blow-up, and $E_x$ the exceptional divisor of $p_x$.
Then one can determine one point $q_x\in E_x\setminus\Delta_0$ on each exceptional divisor whose position depends on $\boldsymbol{\kappa}$ (see \cite{Okamoto}) and blow-up these $4$ points.
We get a smooth projective surface $\Sigma^{\boldsymbol{\kappa}}_t\to\Sigma_0$ inside which the divisor 
$\Gamma=2\Delta+ E_0+ E_1+ E_t+ E_\infty$ is of type $I_0^*$. Then 
$$M_{dR}^{\boldsymbol{\kappa}}(t)\simeq \Sigma^{\boldsymbol{\kappa}}_t\setminus \Gamma.$$
We note that for very special ${\boldsymbol{\kappa}}$, we have to replace $\mathbb P^1\times\mathbb P^1$
by the Hirzebruch surface $\mathbb F^2$ in this construction (see \cite{Okamoto}).

The Riemann-Hilbert correspondence establishes a one-to-one correspondence with the Betti moduli space $M_B^{\boldsymbol{\kappa}}$
of representations
$$\pi_1(\mathbb P^1\setminus\{0,1,t,\infty\})\longrightarrow \mathrm{SL}_2(\mathbb C)$$
by associating to a connection $(E,\nabla)$ in $M_{dR}^{\boldsymbol{\kappa}}(t)$ its monodromy representation.
The Betti moduli space can be constructed as GIT quotient and is represented by an affine surface
$$\mathcal X^{\boldsymbol{\kappa}}=\{(X,Y,Z)\in\mathbb C^3\ ;\ X^2+Y^2+Z^2+XYZ=c_XX+c_YY+c_ZZ+4c\},$$
where $(c_X,c_Y,c_Z,c)\in\mathbb C^4$ depend on $\boldsymbol{\kappa}$ (see \cite{BG,Iwasaki,CL}).
We note that the map $M_B^{\boldsymbol{\kappa}}\to\mathcal X^{\boldsymbol{\kappa}}$
might contract some curves onto singular points for special parameters, namely when $\kappa_x\in\mathbb Z$,
or $\sum_x\pm \kappa_x\in\mathbb Z$; otherwise, $M_B^{\boldsymbol{\kappa}}\simeq\mathcal X^{\boldsymbol{\kappa}}$
is smooth. The Riemann-Hilbert map induces a complex analytic symplectic resolution of the singularities
$$RH^{\boldsymbol{\kappa}}\ :\ M_{dR}^{\boldsymbol{\kappa}}(t)\longrightarrow \mathcal{X}^{\boldsymbol{\kappa}}$$
(see \cite{IIS1,IIS2}). For $\boldsymbol{\kappa}=\left(\frac{1}{2},\frac{1}{2},\frac{1}{2},\frac{1}{2}\right)$, 
we recover the quotient of Atiyah-Serre map $\underline{RH}$ in (\ref{eq:quotientRHAtiyahSerre}).
As in section \ref{sec:Quotient}, the usual compactification $\hat{\mathcal{X}}^{\boldsymbol{\kappa}}$
in $\mathbb P^3$ is smooth, obtained by adding the triangle $XYZ=0$ at infinity; after contracting $Z=0$,
which is a $(-1)$-curve, we get a compactification $\overline{\mathcal{X}}^{\boldsymbol{\kappa}}$
by a half-square type divisor $\underline{D}=\underline{L_1}\cup\underline{L_2}$.

The Riemann-Hilbert map near infinity is studied for instance in \cite{SimpsonAsymptotic,KNPS}.
Let us denote by $W^{\boldsymbol{\kappa}}_t$ and $V^{\boldsymbol{\kappa}}$ the respective neighborhoods 
of $\Gamma\subset \Sigma^{\boldsymbol{\kappa}}_t$ for the deRham side, and $\underline D\subset\overline{\mathcal{X}}^{\boldsymbol{\kappa}}$ in the Betti side.
We propose the following conjecture in the Painlevé VI case:

\begin{conj} The Riemann-Hilbert map near infinity 
$$RH^{\boldsymbol{\kappa}}\ :\ W^{\boldsymbol{\kappa}}_t\setminus \Gamma{\longrightarrow} \underline{V}^{\boldsymbol{\kappa}}\setminus\underline{D}$$
coincides with the classifying map $\underline{LTV}$ of section \ref{sec:TheseGibran}.
\end{conj}

\section{Open questions}

Beyond the conjecture, other questions arise. Denote by $\mathcal M_C$ the moduli space of neighborhoods
of the elliptic curve $C$ with Ueda type $1$, and $\mathcal N$, the moduli space of square-type neighborhoods.

{\bf Generalized isomonodromic deformations.}
For a general neighborhood germ of smooth elliptic curve $(U,C)$, one can deform the classifying map
$$LTV_\tau\ :\ (U,C_\tau)\mapsto (V,D)$$
by deforming the elliptic curve 
$$\mathbb{H}\ni\tau\mapsto C_\tau=\mathbb C/\mathbb Z+\tau\mathbb Z.$$
We derive a classifying map:
$$LTV\ :\ \widetilde{\mathcal{M}}:=\bigcup_{\tau\in\mathbb{H}}\mathcal M_{C_\tau}\to\mathcal N.$$
Taking the quotient by $\Gamma_2\subset\mathrm{SL}_2(\mathbb Z)$, we get a foliation on 
$$\mathcal{M}:=\bigcup_{t\in\mathbb P^1\setminus\{0,1,\infty\}}\mathcal M_{C_t}.$$
The monodromy of this foliation is given by $\Gamma_2$-action on $\mathcal N$ described in \cite[section 1.6]{LTV}
and Painlevé VI neighborhoods give rise to fixed points for this action. It would be nice to understand
if there are other finite orbits. They would give rise to generalized isomonodromic deformations. 

{\bf Other singular elliptic curves.} 
For elliptic curves with automorphisms of order $3$, $4$ or $6$ (not translations),
one might provide similar constructions. For instance, in the case of degree $3$ automorphism,
we will find a map from Atiyah-Serre moduli space (rank $1$ connections) towards
a moduli space of rank $3$ connections with $3$ simple poles on $\mathbb P^1$ for the deRham side.
The compactification, described by Matsumoto in \cite{Matsumoto}, is given by a singular elliptic curve of type $IV$:
three rational curves having self-intersection $(-2)$ and intersecting pairwise transversely at a single point.
On the other hand, the Betti moduli space is described by Lawton in \cite{Lawton} (see also \cite{Komyo}).
It is compactified by a nodal rational curve. But what can be said about other Painlevé cases ?

\end{document}